\documentclass[reqno,11pt]{amsart}
\usepackage[utf8]{inputenc}
\usepackage[margin=1.2 in]{geometry}
\usepackage{amsmath,amssymb,amsthm,mathrsfs,mathtools}
\usepackage{xcolor}
\usepackage{cases}
\usepackage{makecell}
\usepackage{tikz}

\usetikzlibrary{patterns}
\usepackage{cite}
\usepackage[colorlinks=true,citecolor=blue, linkcolor=red,urlcolor=blue]{hyperref}
\usepackage{enumitem}
\numberwithin{equation}{section}
\allowdisplaybreaks[2]

\theoremstyle{plain}
 \newtheorem{theorem}{Theorem}[section]
 \newtheorem{proposition}[theorem]{Proposition}
 \newtheorem{lemma}[theorem]{Lemma}
 
 \newtheorem{corollary}[theorem]{Corollary}
 \theoremstyle{definition}

 \theoremstyle{remark}

\newcommand{\bal}{\begin{align}}
\newcommand{\bbal}{\begin{align*}}
\newcommand{\beq}{\begin{equation}}
\newcommand{\eeq}{\end{equation}}
\newcommand{\bca}{\begin{cases}}
\newcommand{\eca}{\end{cases}}

\newcommand{\pa}{\partial}

\newcommand{\dd}{\mathrm{d}}

\newcommand{\R}{\mathbb{R}}

\numberwithin{equation}{section}
\newcommand{\norm}[2]{\left\lVert #1\right\rVert_{#2}}

\newcommand{\Hinfty}{H^\infty(\R)}
\makeatletter
	\@namedef{subjclassname@2020}{\textup{2020} Mathematics Subject
	Classification}
	\makeatother
\begin{document}
\title[ill-posedness for the generalized Camassa-Holm equation]{Blow-up and norm inflation for the generalized Camassa-Holm equation}

\author{Jinlu Li}
\address{School of Mathematics and Computer Sciences, Gannan Normal University, Ganzhou 341000, China}
\email{lijinlu@gnnu.edu.cn}

\author{Xing Wu}
\address{College of Artificial Intelligence,
Henan Agricultural University, Zhengzhou 450046, China}
\email{ny2008wx@163.com}

\author{Yanghai Yu}
\address{School of Mathematics and Statistics, Anhui Normal University, Wuhu 241002, China}
\email{yuyanghai214@sina.com}

\author{Shuijin Zhang}
\address{School of Mathematics and Computer Sciences, Gannan Normal University, Ganzhou 341000, China}
\email{zhangshuijin@gnnu.edu.cn}

\begin{abstract}
In this paper, we present a simple method to prove that smooth solutions of the generalized Camassa-Holm equation on the real line, associated  with a new class of smooth initial data, not only blows up in finite time but also exhibits  norm inflation in the low-regularity setting. More precisely, by constructing a suitable class of smooth initial data, we establish finite-time blowup of solutions and obtain quantitative upper and lower bounds on the lifespan. Furthermore, we show that for some smooth initial data, the solution exhibits low-regularity $H^s(\mathbb{R})$-norm inflation: the initial data is smooth and arbitrarily small in $H^s(\mathbb{R})$ while the corresponding solution becomes arbitrarily large in $H^s(\mathbb{R})$ after an arbitrarily short time. This proves ill-posedness for the generalized Camassa-Holm equation in low-regularity $H^s(\mathbb{R})$ and improves the results in \cite{Li2023,Li2026}.
\end{abstract}



\date{\today}
\maketitle
\section{Introduction}

In this paper, we consider the Cauchy problem for a generalized Camassa-Holm equation first introduced by Hakkaev and Kirchev \cite{Hakkaev2005}:
\begin{equation}\label{eq1}
\begin{cases}
\displaystyle
u_t-u_{xxt}+\frac{(Q+2)(Q+1)}2u^Qu_x
 =\partial_x\!\left(\frac Q2u^{Q-1}u_x^2+u^Qu_{xx}\right),\\
u(0,x)=u_0(x),\qquad x\in\R,
\end{cases}
\end{equation}
where $Q\ge 1$ is an integer parameter, and $u(t,x)$ denotes the fluid velocity in the spatial direction at time $t\ge 0$.

In the pioneering work of Hakkaev and Kirchev \cite{Hakkaev2005}, they proved local well-posedness in the Sobolev space $H^s(\mathbb{R})$ under the regularity regime $s>\tfrac32$. Subsequent studies extended the well-posedness theory to Besov spaces. More precisely, \cite{Meng2023, Mi2013, Yan2014, 1Yan2014} established Hadamard local well-posedness in $B_{p,r}^s(\mathbb{R})$ for $s>1+1/p$ with $1\le p\le\infty, 1\leq r<\infty$, and also in the critical space $B_{p,1}^{1+1/p}(\mathbb{R})$ with $1\leq p<\infty$. Local ill-posedness was later obtained in \cite{Li2022-2,Meng2023,Deng2024} for the Besov spaces $B_{p,\infty}^s(\mathbb{R})$ with $s>\max\{1+1/p,3/2\}$, $1\leq p \leq \infty$, the critical Besov spaces $B_{p,r}^{1+1/p}(\mathbb{R}),1\leq p, r<\infty$, and $B_{\infty,1}^{1}(\mathbb{R})$.

When $Q=1$, \eqref{eq1} reduces to the celebrated Camassa-Holm (CH) equation:
\begin{eqnarray*}
u_t-u_{xxt}+3uu_x=2u_xu_{xx}+uu_{xxx}.
\end{eqnarray*}
The CH equation first appeared in the work of Fuchssteiner and Fokas \cite{Fokas1981}, and was later rederived as a shallow-water wave equation by Camassa and Holm \cite{Camassa1993}. It exhibits rich geometric and integrable structure: the CH system is completely integrable \cite{Camassa1993, Constantin2001}, admits a bi-Hamiltonian formulation \cite{Fokas1981, Constantin1997}, possesses an infinite family of conserved quantities, and admits explicit peaked travelling profiles $ce^{-|x-ct|}$. Such peakon solutions are well known  as representative profiles for waves of limiting amplitude \cite{Camassa1993, Constantin2007, Constantin2011}. A large body of literature is devoted to the well-posedness and ill-posedness theory for the  Cauchy problem of the CH equation in Sobolev and Besov spaces \cite{Byers2006, Constantin1998, 1Constantin1998, Danchin2001, Danchin2003,Guo2019, Guo2022, Li2016, Li2022, Ye2023}.

A substantial body of work has been devoted to the study of  low-regularity ill-posedness for various Camassa-Holm-type equations. A prominent approach to establishing norm inflation and the failure of continuous dependence on initial data relies on the construction of peakon-antipeakon solutions
$$
u(x, t)=p(t) e^{-|x+q(t)|}-p(t) e^{-|x-q(t)|}.
$$
Byers \cite{Byers2006} utilized this framework to prove norm inflation for the Camassa-Holm equation in $H^s(\mathbb R)$ for $1<s<3/2$. Himonas, Holliman and Grayshan \cite{Himonas2014} established ill-posedness for the Degasperis-Procesi equation in $H^s(\mathbb R)$ for $s<3/2$. For the b-family of equations, Himonas, Grayshan and Holliman \cite{Himonas2016} established ill-posedness via peakon-antipeakon solutions for $b>1$, while Novruzov \cite{Novruzov2021} treated the case $b<1$  under nonvanishing dispersion. Himonas, Holliman and Kenig \cite{Himonas2018} further extended this ill-posedness analysis to the Novikov equation. However, not all Camassa-Holm-type equations admit peakon-antipeakon solutions. For initial data of the form $u_0(x)=p_0 e^{-|x+q_0|}-p_0e^{-|x-q_0|}$, Li, Wu, Yu and Zhu \cite{Li2023} proved that the corresponding solution $u(t,x)\in C([0,T];H^s\cap W^{1,\infty})$ with $s\in [1,\frac32)$ to the Fornberg-Whitham equation exhibits norm inflation. Their analysis relies on the observation that $-\pa_xu(t,x)$ becomes sufficiently large on a nontrivial spatial region. This method was subsequently developed by Li and Yu \cite{Li2026}, who constructed smooth initial data that can be taken arbitrarily small in $H^s$, while the $H^s$ norms of the associated smooth solutions becomes arbitrarily large over an arbitrarily short time interval. We remark that the method in \cite{Li2023,Li2026} applies to the one-dimensional transport equation with  transport velocity $u$:
$$
\pa_tu+uu_x=F(u).
$$
For the one-dimensional transport equation with nonlinear transport velocity $u^Q$:
$$
\pa_tu+u^Qu_x=F(u),
$$
the initial data constructed in \cite{Li2023,Li2026} cannot guarantee that $-\pa_xu(t,x)$ becomes sufficiently large on a nontrivial spatial region.

We first rewrite \eqref{eq1} into the following nonlocal form
\begin{equation}\label{gch}
\begin{cases}
\displaystyle
u_t+u^Qu_x=-G_x*P_Q(u),\\
u(0,x)=u_0(x),\qquad x\in\R,
\end{cases}
\end{equation}
where $G(x)=\tfrac12e^{-|x|}$ is the kernel of inverse operator $(1-\partial_x^2)^{-1}$ and $*$ denotes convolution, and we dentoe
\begin{align*}
A_Q=\frac{Q(Q+3)}{2(Q+1)},\qquad
P_Q(u)=A_Qu^{Q+1}+\frac Q2u^{Q-1}u_x^2.
\end{align*}

In this paper, we investigate the Cauchy problem for the generalized Camassa-Holm equation \eqref{gch}.
Now we state the main results of this paper.

\begin{theorem}[Finite-time blow-up for smooth initial data]\label{th1} Let $Q\geq1$, $1<s<3/2$ and let $u_{0,n}$ be defined by \eqref{id}.
Then the solution to the generalized Camassa-Holm equation \eqref{gch} blows up in finite time. Moreover, let $T^*_n$ denote the lifespan of the
solution to \eqref{gch} with initial data $u_{0,n}$. Then
\bbal
 T_n^-\leq T_n^*\leq T_n^+,
\end{align*}
where $$T_n^-\approx T_n^+\approx\frac2{Q}n^{\kappa(Q-1)-1},\quad 0<\kappa<\min\left\{\frac{1}{3-2s}-1,\frac1Q\right\}.$$
\end{theorem}
\begin{theorem}[Norm inflation for smooth initial data]\label{th2}
Let $Q\geq1$ and $1<s<3/2$ be fixed.
For any $\varepsilon>0$, there exists real-valued initial data
$u_0\in C_c^\infty(\R)\subset\mathcal S(\R)$ and a time
$0<T_0<\varepsilon$
such that the solution $u\in C([0,T_0];\Hinfty)$ to the generalized Camassa-Holm equation \eqref{gch} satisfies
\begin{align*}
 \norm{u_0}{H^s}\leq\varepsilon\qquad\text{but}
 \qquad \norm{u(T_0)}{H^s}>\varepsilon^{-1}.
\end{align*}
Therefore, smooth solutions exhibit norm inflation near the zero initial data in $H^s(\R)$ whenever $1<s<3/2$.
\end{theorem}

Following the proof of Theorem \ref{th2} and suitably modifying the initial data, we can show that smooth solutions to \eqref{gch} exhibit norm inflation near the zero initial data in either $B^s_{p,r}(\mathbb{R})$ or $W^{1,p}(\mathbb{R})$.

\begin{corollary}
Let $1<s<1+\frac1p$ and $(p,r)\in [2,\infty)\times [1,\infty]$.
For any $\varepsilon>0$, there exists a solution $u\in \mathcal{C}([0,T_0];B^s_{p,r}(\mathbb{R}))$ to \eqref{gch} with initial datum $u_0$ such that
$$\|u_0\|_{B^s_{p,r}}\leq \varepsilon\qquad\text{but}
 \qquad  \|u(T_0)\|_{B^s_{p,r}}\geq\varepsilon^{-1},$$
for some $0<T_0<\varepsilon$.
\end{corollary}
\noindent{\bf Sketch of the Key Ideas.}\; As mentioned above,
for the one-dimensional transport equation with nonlinear transport velocity $u^Q$:
$$
\pa_tu+u^Qu_x=F(u),
$$
the initial data constructed in \cite{Li2023,Li2026} cannot guarantee that $-\pa_xu(t,x)$ becomes sufficiently large on a nontrivial spatial region. For initial data similar to those in \cite{Li2023,Li2026}, the quantity $u^Q$ turns out to be zero, which suppresses the growth of $-\pa_xu(t,x)$. To overcome this difficulty, we first modify the high-frequency profile $g_n(x)=nq_n\psi(x/q_n)$ to replace the initial data used in \cite{Li2023,Li2026}. One can verify that $g_n$ satisfies
$$
\partial_x g_n(x)=-n \quad \text{on } (-q_n,q_n),
$$
and
$$\|u_{0,n}\|_{H^s} \simeq \|\partial_x u_{0,n}\|_{L^{p_s}(-q_n,q_n)} \quad \text{with} \quad p_s=\dfrac{1}{3/2-s}.
$$
Next, we introduce a low-frequency function $f_n(x)=\delta_n\chi(x)$ to ensure that the solution $u(t,x)$ remains close to $\delta_n$ for $t\leq T^*$. For this initial data $u_{0,n}=f_n+g_n$, we prove that $-\pa_xu(t,x)$ becomes sufficiently large on a nontrivial spatial region. Finally, following the strategy developed in \cite{Li2023,Li2026}, we show that the norm inflation effect is controlled by the quantity
$$
\left(\int_{\phi_n(t,-q_n)}^{\phi_n(t,q_n)} |\partial_x u_n|^{p_s}\,\mathrm{d}x\right)^{\frac{1}{p_s}}.
$$
\noindent{\bf Outline of the paper.} Section \ref{sec2} gathers notation and recalls auxiliary results needed throughout the paper. Sections \ref{sec3} and \ref{sec4} are devoted to the proof of the main theorems, following the steps outlined below:
\begin{enumerate}
\item Construction of smooth initial data;
\item Estimation for $\frac Q2 u_n^{Q-1}(t,x)$;
\item Upper and lower bounds on the lifespan;
\item Pre-blow-up time;
\item Moving-interval $L^p$ norm of $-\pa_xu_n$;
\item Norm inflation.
\end{enumerate}

\section{Notations and preliminaries}\label{sec2}
\subsection{Notation}
Given a Banach space X, we denote its norm by $\|\cdot\|_{X}$. For $I\subset\mathbb{R}$, let $\mathcal{C}(I;X)$ stand for the space of continuous X-valued functions on $I$. For positive sequences $\{a_n\}$ and $\{b_n\}$, we write $a_n\approx b_n$ to mean that
$$a_n = b_n + o(b_n),\qquad \lim_{n\to\infty}\frac{o(b_n)}{b_n}=0.$$
For every $f\in\mathcal{S}'(\mathbb{R})$, the Fourier transform $\widehat{f}$ on $\mathbb{R}$ is defined as
$$\widehat{f}(\xi)=\int_{\mathbb{R}}e^{-ix\xi}f(x)\,\dd x,\quad \forall\xi\in\mathbb{R}.$$
For $s\in\mathbb{R}$, the Bessel potential operator $J^{s}:=(1-\partial_x^2)^{\frac{s}{2}}$ is specified via its Fourier multiplier:
$$\widehat{J^s f}(\xi)=(1+|\xi|^2)^{\frac{s}{2}}\widehat{f}(\xi).$$
For $s\in\mathbb{R}$, the inhomogeneous Sobolev space $H^s(\mathbb{R})$ is equipped with the norm
$$\|f\|_{H^s(\mathbb{R})}^2=\|J^{s}f\|_{L^2(\R)}^2=\int_{\mathbb{R}}(1+|\xi|^2)^s|\widehat{f}(\xi)|^2\,\dd\xi.$$
For $f\in H^s(\mathbb{R})$ with $s\ge 0$, we define the homogeneous Sobolev semi-norm on $\dot{H}^s(\mathbb{R})$ by
$$\|f\|_{\dot{H}^s(\mathbb{R})}^2=\int_{\mathbb{R}}|\xi|^{2s}|\widehat{f}(\xi)|^2\,\dd\xi.$$
Next, we recall several known results concerning solutions to the generalized Camassa-Holm equation \eqref{gch}, including
local well-posedness, $H^1$ conservation law and $L^\infty$ boundedness.
\subsection{Local Well-posedness}
First, we recall the local well-posedness result for the generalized Camassa-Holm equation \eqref{gch} in Sobolev spaces, see \cite{Hakkaev2005}.
\begin{lemma}\label{le1} Let $u_0 \in H^\sigma(\mathbb{R})$ with $\sigma>\frac{3}{2}$. Then there exists  $T^*>0$ such that the generalized Camassa-Holm equation \eqref{gch} admits a unique solution
$$
u\in C\left([0, T^*) ; H^\sigma(\mathbb{R})\right) \cap C^1\left([0, T^*); H^{\sigma-1}(\mathbb{R})\right).
$$
\end{lemma}
The maximum value of $T^*$ in Lemma \ref{le1} is called the maximum existence time or the lifespan of the solution in general. If $T^*<\infty$, that is $\limsup_{t \rightarrow T^*} \|u(t)\|_{H^\sigma}=\infty$, we say that the solution blows up in finite time.

\subsection{$H^1$ conservation law, blow-up criterion and $L^\infty$ boundedness}

\begin{lemma}\label{lem:energy}
Assume that  $u\in C([0,T^*);H^{\sigma})$ is a solution to the generalized Camassa-Holm equation \eqref{gch}. Then for all $t\in [0,T^*),$ we have
\begin{equation}\label{gch-H1L}
 \norm{u(t)}{H^1}= \norm{u_0}{H^1}=:E_0,
 \qquad \norm{u(t)}{L^\infty}\leq E_0/\sqrt2.
\end{equation}
\end{lemma}
\begin{proof}
Multiplying \eqref{eq1} by $u$ and integrating over $\R$, we obtain
\begin{align*}
 \frac12\frac{\mathrm d}{\mathrm dt}
      \int_\R(u^2+u_x^2)\dd x
 &=-\frac Q2\int_\R u^{Q-1}u_x^3\dd x
   -\int_\R u^Qu_xu_{xx}\dd x\\
 &=-\frac Q2\int_\R u^{Q-1}u_x^3\dd x
   +\frac Q2\int_\R u^{Q-1}u_x^3\dd x=0,
\end{align*}
which implies the first part of \eqref{gch-H1L}. Using the one-dimensional inequality
$\norm{f}{L^\infty}^2\leq\tfrac12\norm{f}{H^1}^2$, we obtain the second estimate in \eqref{gch-H1L}.
This completes the proof of the lemma.
\end{proof}
Next, we establish a blow-up criterion for \eqref{gch} and derive estimates for \(F_Q(u)\) and \(R_Q(u)\) by exploiting  the conservation of the \(H^1\) energy.

\begin{lemma}\label{lem:blowup}
Let $u_0\in H^\sigma(\mathbb{R})$ with $\sigma>\frac{3}{2}$, and let $T^*$ be the maximal existence time of the
solution $u$ to \eqref{gch} with the initial data $u_0$. Then the corresponding solution blows up in finite time if and only
\begin{align*}
  \limsup_{t \rightarrow T^*} \|\partial_x u(t)\|_{L^\infty}=\infty.
\end{align*}
\end{lemma}

\begin{proof}
On the one hand, if \(\limsup_{t \rightarrow T^*} \|\partial_x u\|_{L^\infty}=\infty\), then \(\limsup_{t \rightarrow T^*} \|u\|_{H^\sigma}=\infty\). On the other hand, assuming that \(\limsup_{t \rightarrow T^*} \|u\|_{H^\sigma}=\infty\), we  prove that \(\limsup_{t \rightarrow T^*} \|\partial_x u\|_{L^\infty}=\infty\) by contradiction. Suppose, on the contrary, that
\(\limsup_{t \rightarrow T^*} \|\partial_x u\|_{L^\infty}<\infty\). First, using standard energy estimates, we obtain
\begin{align*}
\frac{d}{dt}\|u\|^2_{H^\sigma}
&\leq C\|u\|_{L^\infty}^{Q-1}\|\partial_x u\|_{L^\infty}\|u\|^2_{H^\sigma}\\
&\leq C\|u_0\|_{H^1}^{Q-1}\|\partial_x u\|_{L^\infty}\|u\|^2_{H^\sigma}.
\end{align*}
Here we have used the Sobolev embedding $H^1(\mathbb{R})\hookrightarrow L^\infty(\mathbb{R})$ and  conservation of the $H^1$ energy  for $u$ in the second inequality.

By Gronwall’s inequality, we obtain
\begin{align*}
\|u\|^2_{H^\sigma}
&\leq \|u_0\|_{H^\sigma}^2\exp\left(C\|u_0\|_{H^1}^{Q-1}\int_0^{T^*}\|\partial_x u\|_{L^\infty}\,dt\right)\\
&\leq \|u_0\|_{H^\sigma}^2\exp\left(C T^*\|u_0\|_{H^1}^{Q-1}\sup_{t\in[0,T^*)}\|\partial_x u(t)\|_{L^\infty}\right).
\end{align*}
As \(\limsup_{t \rightarrow T^*}\|\partial_x u\|_{L^\infty}<\infty\), the exponential factor remains finite, which implies \(u\in L^\infty([0,T^*];H^\sigma)\). This contradicts the maximality of \(T^*\).
This completes the proof of the lemma.
\end{proof}

\begin{lemma}\label{lem:energy2}
Define
\begin{align*}
 F_Q(u)=-G_x*P_Q(u),\qquad
 R_Q(u)=A_Qu^{Q+1}-G*P_Q(u).
\end{align*}
Then
\begin{equation}\label{es-FQRQ}
 \max\{\norm{F_Q(u(t))}{L^\infty},
        \norm{R_Q(u(t))}{L^\infty}\}\leq K:=A_QE_0^{Q+1}.
\end{equation}
\end{lemma}
\begin{proof}
Set $M=\norm{u(t)}{L^\infty}$. Since $A_Q\geq Q/2$, we deduce
\begin{align*}
 \norm{P_Q(u)}{L^1}
 \leq A_QM^{Q-1}\bigl(\norm{u}{L^2}^2+\norm{u_x}{L^2}^2\bigr)
 =A_QM^{Q-1}E_0^2,
\end{align*}
which yields the first part of \eqref{es-FQRQ}. Using $\norm{G}{L^\infty}=\norm{G_x}{L^\infty}=1/2$ and $M^2\leq E_0^2/2$, we obtain
\begin{align*}
 \norm{F_Q(u)}{L^\infty}
 \leq\tfrac{A_Q}{2}M^{Q-1}E_0^2\leq A_QE_0^{Q+1},
\end{align*}
and
\begin{align*}
 \norm{R_Q(u)}{L^\infty}
 &\leq A_QM^{Q+1}+\tfrac{A_Q}{2}M^{Q-1}E_0^2
 \leq A_QM^{Q-1}E_0^2\leq A_QE_0^{Q+1}.
\end{align*}
This completes the proof of the lemma.
\end{proof}

\subsection{The Equation Along the Flow}
Given a Lipschitz velocity field $u$, we may solve the following ODE to find the flow map induced by $u$:
\begin{align}\label{phi-flow}
\quad\begin{cases}
 \partial_t\phi(t,x)=u^Q(t,\phi(t,x)),\\
\phi(0,x)=x.
\end{cases}
\end{align}
Since the velocity field is Lipschitz,  we have  for all $(t, x) \in[0,T^*)\times\R$
\bbal
\pa_x\phi(t,x)=\exp\left(\int_0^t Q u^{Q-1}\pa_xu(\tau,\phi(\tau,x))\dd\tau\right)>0.
\end{align*}
This shows that $\phi(t,\cdot)$ is an increasing diffeomorphism on $\R$. That is,  for all $x,y\in \R,$  $\phi(t,x)<\phi(t,y)$ whenever $x< y$. Applying a classical result from the theory of ordinary differential equations, we obtain

\begin{lemma}\label{le4} Assume $u_0 \in H^\sigma(\mathbb{R})$ with $\sigma>\frac{3}{2}$. If $T^*$ is defined as in Lemma \ref{le1}, then the ODE \eqref{phi-flow} has a unique solution $\phi \in C^1([0, T^*) \times \mathbb{R})$. Moreover, for every $t\in[0,T^*)$, the map $\phi(t, \cdot):\R\rightarrow\R$ is an increasing diffeomorphism on $\mathbb{R}$ and $\pa_x\phi(t, x)>0$ on $[0, T^*) \times \mathbb{R}$.
\end{lemma}
\begin{lemma}[Characteristic equations]\label{le-UW}
Let $\phi(t,x)$ denote the flow defined by \eqref{phi-flow}.
Setting $U(t,x)=u(t,\phi(t,x))$ and
$W(t,x)=(\pa_xu)(t,\phi(t,x))$, we have
\begin{equation}\label{eq-U}
\begin{cases}
 \partial_t U=F_Q(u)(t,\phi),\\
 U(t=0,x)=u_0(x),
\end{cases}
\end{equation}
and
\begin{equation}\label{eq-W}
\begin{cases}
\partial_t W=-\frac Q2U^{Q-1}W^2+R_Q(u)(t,\phi),\\
  W(t=0,x)=u'_0(x).
\end{cases}
\end{equation}
In particular, for $t\in[0,T^*)$, we have
\begin{equation}\label{UleqK}
 \|U(t,x)-u_0(x)\|_{L^\infty(\R)}\leq Kt.
\end{equation}
\end{lemma}

\begin{proof}
Differentiating \eqref{gch} and using
$-\partial_x^2(G*f)=f-G*f$, we derive
\begin{equation}\label{eq-u-ux}
\begin{cases}
u_t+u^Qu_x=-G_x*P_Q(u),\\
u_{tx}+u^Qu_{xx}=-\frac Q2u^{Q-1}u_x^2+R_Q(u).
\end{cases}
\end{equation}
Combining \eqref{eq-u-ux} with the flow map \eqref{phi-flow}, we arrive at \eqref{eq-U}-\eqref{eq-W}. Integrating \eqref{eq-U} and applying Lemma \ref{lem:energy}, we obtain \eqref{UleqK}.
\end{proof}

\section{Proof of Theorem \ref{th1}: Finite-time blow-up}\label{sec3}

{\bf Choice of Parameters}. From now on, we fix $Q\geq1$ and $1<s<3/2$.  Throughout this paper, the parameters are chosen as follows:
\begin{align*}
& p_s=\frac1{3/2-s}>2,\\
&q_n=(n\log n)^{-p_s},\\
&0<\kappa<\min\!\left\{\frac{p_s-2}{2},\frac1Q\right\}, \\
&\delta_n=n^{-\kappa},\quad n\gg1,\\
&\rho_n=n^{-1}\delta_n=n^{-\kappa-1}.
\end{align*}

\subsection{The Example Class for Smooth Initial Data}
First, we construct an example of smooth initial data as follows.
Choose an even function $\chi\in C_c^\infty((-2,2))$ satisfying
$0\leq\chi\leq1$ and $\chi=1$ on $[-1,1]$. Set
\begin{align*}
 &h(x)=\tfrac12\chi(x-5)+\tfrac12\chi(x+5)-\chi(x),\\
&\psi(x)=\int_0^x h(y)\dd y.
\end{align*}
The supports of the above three terms defining $h$ are disjoint. Thus
$h$ is even, $|h|\leq1$, $h=-1$ on $[-1,1]$, and
$\int_\R h(x)\dd x=0$. Moreover,
\begin{align*}
 \psi\in C_c^\infty((-7,7)),\qquad
 \psi(x)=-x\quad(|x|\leq1),\qquad
 \norm{\psi'}{L^\infty(\R)}=1.
\end{align*}
Define the initial data as the sum of a high-frequency part $g_n=nq_n\psi(x/q_n)$ and a low-frequency part $f_n=\delta_n\chi(x)$:
\begin{equation}\label{id}
 u_{0,n}(x)=f_n+g_n=\delta_n\chi(x)+nq_n\psi\left(\frac{x}{q_n}\right).
\end{equation}
For all sufficiently large $n$, whenever $|x|\leq q_n$,
\begin{equation}\label{id-qn}
 u_{0,n}(x)=\delta_n-nx,\qquad
 u_{0,n}'(x)=-n.
\end{equation}
It also follows from \eqref{id}--\eqref{id-qn} that $\norm{u_{0,n}'}{L^\infty(\R)}=n$ since $\delta_n\norm{\chi'}{L^\infty}<n$.

Furthermore, we have the following lemma.
\begin{lemma}\label{le001} Let $u_{0,n}$ be defined by \eqref{id} and $n\gg 1$. For $r\geq0$, there exist constants $C,C_E>0$ independent of $n$ such that
\begin{align}\label{id-Hrs1}
\|u_{0,n}\|_{H^r}\leq C (nq_n^{3/2-r}+\delta_n).
\end{align}
In particular,
\begin{align*}
\|u_{0,n}\|_{H^s}\leq \frac{C }{\log n}.
\end{align*}
and
\begin{align*}
\|u_{0,n}\|_{H^1}\leq C_E \delta_n=C_En^{-\kappa}.
\end{align*}
\end{lemma}
\begin{proof} By Fourier scaling, it is straightforward to verify that the high-frequency function $g_n$ satisfies
\begin{align*}
\|g_n\|_{H^r}^2&=n^2q_n^{4}\int_\R(1+|\xi|^2)^r|\widehat\psi(q_n\xi)|^2\dd\xi
=n^2q_n^{3-2r}\int_\R(q_n^2+|\eta|^2)^r|\widehat\psi(\eta)|^2\dd\eta.
\end{align*}
Hence, for $n\gg 1$,
$$
n^2q_n^{3-2r}\|\psi\|_{\dot{H}^r}^2 \leq \|g_n\|_{H^r}^2\leq n^2q_n^{3-2r}\|\psi\|_{H^r}^2,
$$
which implies the first estimate in \eqref{id-Hrs1}. In particular, for $n\gg 1$,
\begin{align*}
 \norm{u_{0,n}}{H^s} \leq \delta_n \norm{\chi}{H^s}+nq_n^{3/2-s}\norm{\psi}{H^s}
 \leq \left(n^{-\kappa}+\frac1{\log n}\right)\max\{\norm{\chi}{H^s},\,\norm{\psi}{H^s}\},
\end{align*}
and
\begin{align*}
 \norm{u_{0,n}}{H^1}
 \leq \delta_n\norm{\chi}{H^1}+n^{1-\frac12p_s+\kappa}(\log n)^{-\frac12p_s} \delta_n\norm{\psi}{H^1}
 \leq \delta_n\max\{\norm{\chi}{H^1},\,\norm{\psi}{H^1}\},
\end{align*}
where we have used the restriction $\kappa<\frac12p_s-1$. This completes the proof of the lemma.
\end{proof}

\subsection{Estimation of $\frac Q2 u^{Q-1}_n(t,x)$}
From now on, we assume that the Cauchy problem \eqref{gch} with initial data $u_{0,n}$ given by \eqref{id}
 has a unique solution $u_n\in C([0, T_n^*); H^\infty(\R))$, where $T^*$ is the lifespan of the solution and $H^\infty(\R)=\cap_{s=1}^\infty H^s(\R)$. In fact, one may combine  a priori estimate and a compactness argument to  establish local-in-time well-posedness in $H^2(\R)$. From \eqref{gch} and standard product estimates in Sobolev spaces, we deduce that $u_n\in C^1([0, T_n^*); H^\infty(\R))$.
This suffices to render the subsequent arguments rigorous.
Next, we show that $\frac Q2u^{Q-1}_n(t,x)\approx \frac Q2\delta_n^{Q-1}$ for $t\leq \frac4{Qn^{1-\kappa(Q-1)}}$ and $x\in[\phi_n(t,-q_n),\phi_n(t,q_n)]$.
First, we need to establish the following estimate for $u_n$.
\begin{lemma}\label{es-un}
Let $n\gg 1$. There exists a constant $C_1>0$ independent of $n$ such that for $t< \frac4{Qn^{1-\kappa(Q-1)}}$,
\begin{align*}
\|u_n(t,\cdot)\|_{L^\infty(\R)}\leq \delta_n(1+C_1\rho_n),
\end{align*}
and
\begin{align*}
u_n(t,\phi_n(t,x))&\geq  \delta_n(1-C_1\rho_n),\quad |x|\leq q_n.
\end{align*}
\end{lemma}
\begin{proof}
Let $\phi_n$ be the flow associated with $u_n^Q$.
Fix $C_0= A_QC_E^{Q+1}$ and define
\begin{align*}
 &K_n=C_0\delta_n^{Q+1}=C_0n^{-\kappa(Q+1)},\\
 &a_n=\frac Q2\delta_n^{Q-1}=\frac Q2n^{-\kappa(Q-1)},\\
 &\overline T_n=\frac2{a_nn}=\frac4{Qn^{1-\kappa(Q-1)}}.
\end{align*}
By Lemma \ref{lem:energy}, one has
\begin{align*}
 \max\left\{\|F_Q(u_n(t))\|_{L^\infty},
        \|R_Q(u_n(t))\|_{L^\infty}\right\}\leq K_n.
\end{align*}
For $0\leq t<\min\{T_n^*,\overline T_n\}$, Lemma \ref{le-UW} gives
\begin{equation}\label{un-u0n}
 |u_n(t,\phi_n(t,x))-u_{0,n}(x)|
 \leq K_nt\leq\frac{4C_0}{Q}\delta_n\rho_n.
\end{equation}
Since $\kappa<\frac{p_s-2}{2}$, we have $n^{1-p_s}(\log n)^{-p_s}=nq_n\leq \delta_n\rho_n=n^{-2\kappa-1}$ for large $n$. Set $C_1=\frac{4C_0}{Q}+\|\psi\|_{L^\infty}>0$. Using \eqref{un-u0n}, for all $x\in \R$,
\begin{align*}
 |u_n(t,\phi_n(t,x))|
 &\leq |u_{0,n}(x)|+\frac{4C_0}{Q}\delta_n\rho_n\nonumber\\
 &\leq \delta_n+nq_n\|\psi\|_{L^\infty}+\frac{4C_0}{Q}\delta_n\rho_n\nonumber\\
 &\leq \delta_n(1+C_1\rho_n),
\end{align*}
and for $|x|\leq q_n$,  $n\gg 1$,
\begin{align*}
 u_n(t,\phi_n(t,x))&\geq u_{0,n}(x)-\frac{4C_0}{Q}\delta_n\rho_n \\&
 \geq  \delta_n-nq_n\|\psi\|_{L^\infty}-\frac{4C_0}{Q}\delta_n\rho_n \\&
 \geq  \delta_n(1-C_1\rho_n).
\end{align*}
This completes the proof of the lemma.
\end{proof}

As a direct consequence for $\frac Q2u^{Q-1}_n$, we have
\begin{corollary}\label{co-unQ-1}
Let $n\gg 1$. There exists a constant $C_2>0$ independent of $n$ such that for $t< \frac4{Qn^{1-\kappa(Q-1)}}$
\begin{align*}
\left\|\frac Q2u^{Q-1}_n(t,\cdot)\right\|_{L^\infty(\R)}\leq a_n^+:=a_n(1+C_2\rho_n),
\end{align*}
and for $|x|\leq q_n$
\begin{align*}
\frac Q2 u^{Q-1}_n(t,\phi_n(t,x))\geq a_n^-:= a_n(1-C_2\rho_n).
\end{align*}
\end{corollary}
\begin{proof}
Set $C_2=2(Q-1)C_1$. Then, by Lemma \ref{es-un}, for $x\in \R$, we have
\begin{align*}
 \left|\frac Q2u^{Q-1}_n(t,x)\right|&\leq \frac Q2\delta_n^{Q-1}(1+C_1\rho_n)^{Q-1}
 \\&\leq a_n(1+C_1\rho_n)^{Q-1}
 \\ &\leq a_n\left(1+2(Q-1)C_1\rho_n\right)\leq a_n^+,
\end{align*}
and for $|x|\leq q_n$
\begin{align*}
\frac Q2u^{Q-1}_n(t,\phi_n(t,x))&\geq \frac Q2\delta_n^{Q-1}(1-C_1\rho_n)^{Q-1} \\&
\geq a_n(1-C_1\rho_n)^{Q-1} \\&
\geq a_n\left(1-\frac12(Q-1)C_1\rho_n\right)\geq a_n^-,
\end{align*}
where we have used the simple facts that $(1+y)^{Q-1}\leq 1+2(Q-1) y$ and $ (1-y)^{Q-1}\geq 1-\frac12(Q-1) y$  for small $y\geq 0$.
\end{proof}

\subsection{Upper and Lower Bound of Lifespan}
We now estimate the maximal existence time $T^*_n$ of the solution $u_n$ with initial data $u_{0,n}$.
\begin{proposition}\label{le-ul-T}
Let $n\gg 1$. The lifespan $T^*_n$ of the solution $u_n$ to the generalized Camassa-Holm equation \eqref{gch} with initial data $u_{0,n}$ given by \eqref{id} satisfies
$$
 \frac1{a_n^+(n+d_n^+)}:= T_n^-\leq T_n^*\leq T_n^+:= \frac1{a_n^-(n-d_n^-)},
$$
where
\begin{align}
 &d_n^+=\sqrt{K_n/a_n^+}=\sqrt{\frac{2C_0}{Q}} n^{-\kappa}+o(n^{-\kappa}), \label{dn+}\\
 &d_n^-=\sqrt{K_n/a_n^-}=\sqrt{\frac{2C_0}{Q}} n^{-\kappa}+o(n^{-\kappa}).\label{dn-}
\end{align}
\end{proposition}
\begin{proof}\,
{\bf \underline{Upper Bound of Lifespan}.}
Now we define
$$
X_n(t,x)=-(\pa_xu_{n})(t,\phi_n(t,x)),\quad |x|\leq q_n.
$$
Then $X_n(0,x)=n$ for $|x|\leq q_n$. By Lemma \ref{le-UW} and Corollary \ref{co-unQ-1}, from \eqref{eq-W}, we obatin
\begin{equation}\label{zn:geq}
 \partial_tX_n\geq a_n^-X_n^2-K_n,\quad \forall (t,x)\in[0,\min\{T_n^*,\overline T_n\})\times [-q_n,q_n].
\end{equation}
Since $X_n(0,x)=n>d_n^-=\sqrt{K_n/a_n^-}\approx \sqrt{\frac{2C_0}{Q}}n^{-\kappa}$ for $|x|\leq q_n$ when $n$ is large, a continuity argument yields
$$
X_n(t,x)\geq X_n(0,x)=n, \quad \forall (t,x)\in[0,\min\{T_n^*,\overline T_n\})\times [-q_n,q_n].
$$
Note that $K_n=a_n^-(d_n^-)^2$, we rewrite
\begin{align}\label{zn(t)-di}
 \partial_tX_n\geq a_n^-(X_n-d_n^-)^2, \quad \forall (t,x)\in[0,\min\{T_n^*,\overline T_n\})\times [-q_n,q_n].
\end{align}
Solving this differential inequality yields
\begin{equation}\label{zn:lower}
 X_n(t,x)\geq d_n^-+
       \frac{n-d_n^-}{1-a_n^-(n-d_n^-)t},
\end{equation}
for $t<\min\{T_n^*,\overline T_n\}$. Since $C_1\rho_n\to0$ and $d_n^-/n\to0$, we have $\frac1{na_n} \approx T_n^+<\overline T_n=\frac2{na_n}$ for $n\gg 1$.
If $T_n^*>T_n^+$, then \eqref{zn:lower} contradicts the regularity at
$T_n^+$.
Thus the solution to \eqref{gch} must blow up in finite time and $T_n^*\leq T_n^+$.

{\bf \underline{Lower Bound of Lifespan}.}
By \eqref{dn+}-\eqref{dn-}, there exists a fixed $C_4=2\sqrt{\frac{2C_0}{Q}}$ such that for large $n$
\begin{equation}\label{eq:d-estimate}
 0<d_n^\pm\leq C_4\delta_n= C_4n^{-\kappa}\;\Rightarrow\;
 \frac{d_n^\pm}{n}\leq C_4\rho_n.
\end{equation}
Denote $W_n(t,x)=(\pa_xu_n)(t,\phi_n)$. By Lemma \ref{le-UW} and Corollary \ref{co-unQ-1}, from \eqref{eq-W}, we have
\begin{align}\label{Wn(t)}
 \frac{\mathrm d}{\mathrm dt}W_n(t,x)
 = -\frac Q2 u^{Q-1}_n(t,\phi_n)W^2_n(t,x)+R_Q(u_n)(t,\phi_n),\quad W_n(0,x)=u'_{0,n}(x),
\end{align}
with
$$
\left|\frac Q2 u^{Q-1}_n(t,\phi_n)\right|\leq a_n^+, \qquad |R_Q(u_n)(t,\phi_n)|\leq K_n.
$$
Introduce the key comparable quantity
\begin{equation}\label{Yn(t)}
 Y_n(t)=\frac{n+d_n^+}{1-a_n^+(n+d_n^+)t}-d_n^+\quad\text{with}\quad Y_n(0)=n.
\end{equation}
Obviously,
\begin{align*}
 Y_n'=a_n^+(Y_n+d_n^+)^2\geq a_n^+Y_n^2+K_n.
\end{align*}
From this and \eqref{Wn(t)}, one has for $t<T_n^*$
\begin{align*}
(Y_n-W_n)'&=a_n^+(Y_n+d_n^+)^2+\frac Q2U^{Q-1}W_n^2-R_Q(u_n)(t,\phi)\\&
\geq a_n^+Y_n^2+K_n-a_n^+W_n^2-K_n \\&
\geq a_n^+(Y_n-W_n)(Y_n+W_n).
\end{align*}
A similar argument yields that for $t<T_n^*$
$$
(Y_n+W_n)'\geq a_n^+(Y_n-W_n)(Y_n+W_n).
$$
Since $Y_n(0)\pm W_n(0)\geq 0$, we obtain  for $t<T_n^*$
$$
(Y_n-W_n)\exp\left(-\int^t_0 a_n^+(Y_n+W_n)\,\mathrm d \tau\right)\geq Y_n(0)- W_n(0)\geq 0,
$$
and
$$
(Y_n+W_n)\exp\left(-\int^t_0 a_n^+(Y_n-W_n)\,\mathrm d \tau\right)\geq Y_n(0)+ W_n(0)\geq 0,
$$
which implies $|W_n|\leq Y_n$. Hence for $t<T_n^*$
\begin{align}\label{ux-Yn}
 \norm{\pa_xu_{n}(t,\cdot)}{L^\infty}\leq Y_n(t).
\end{align}
If $T_n^*<T_n^-$, the right-hand side of \eqref{ux-Yn} remains bounded up to $T_n^*$, which yields a contradiction. Therefore, we conclude that $T_n^*\geq T_n^-$.

Combining the two cases above, we complete the proof of Proposition \ref{le-ul-T} and hence establish Theorem \ref{th1}.
\end{proof}
\section{Proof of Theorem \ref{th2}: Norm Inflation}\label{sec4}

\subsection{Pre-blow-up time}
We now choose $t_n<T_n^-$ approaching $T_n^-$ such that $(1-r_nt_n)^{-1}$ becomes large and yields the growth of the $H^s$ norm of $u_n(t,x)$.
\begin{lemma}\label{le-pre-T}
There exists a time $t_n<T_n^-$ and a constant $c_6\in(0,1)$, independent of $n\gg 1$, such that
$$
\frac{c_6}{\rho_n}\leq (1-r_nt_n)^{-1}\leq \frac{1}{\rho_n},\qquad
 r_n:=a_n^-(n-d_n^-),.
$$
\end{lemma}
\begin{proof}
Recall that $T_n^-$, we compute
\begin{align}\label{1-rnTn}
1-r_nT_n^-&=1- \frac{a_n^-(n-d_n^-)}{a_n^+(n+d_n^+)} \nonumber
\\&=1-\frac{1-C_2\rho_n}{1+C_2\rho_n}\cdot\frac{1-\frac{d_n^-}{n}}{1+\frac{d_n^+}{n}} \nonumber
\\&=\frac{(1+C_2\rho_n)(1+\frac{d_n^+}{n})-(1-C_2\rho_n)(1-\frac{d_n^-}{n})}{(1+C_2\rho_n)(1+\frac{d_n^+}{n})} \nonumber
\\&=\frac{2C_2\rho_n+\frac{d_n^+}{n}+\frac{d_n^-}{n}+C_2\rho_n\frac{d_n^+-d_n^-}{n}}{(1+C_2\rho_n)(1+\frac{d_n^+}{n})},
\end{align}
which is comparable to $(2C_2+\frac{4C_0}{Q})\rho_n$ since $\frac{d_n^+}{n}  \approx \frac{d_n^+}{n} \approx  \sqrt{\frac{2C_0}{Q}} n^{-\kappa-1}=\sqrt{\frac{2C_0}{Q}}\rho_n$. Now we set
\begin{equation}\label{pro-time}
 t_n=(1-\rho_n)T_n^-.
\end{equation}
Since $\kappa<1/Q$, then
$$
 t_n \approx \frac1{a_nn}
 =\frac2Qn^{-1+\kappa(Q-1)}\longrightarrow0,
$$
and
$$
1-r_nt_n=1-r_n(1-\rho_n)T_n^-=\rho_n+(1-\rho_n)
(1-r_nT_n^-).
$$
From \eqref{eq:d-estimate} and \eqref{1-rnTn}, one has
$$
0\leq1-r_nT_n^-\leq \frac{2C_2\rho_n+2d_n^+/n+C_2\rho_nd_n^+/n}{(1+C_2\rho_n)(1+d_n^+/n)}\leq(2C_2+2C_4)\rho_n.
$$
Hence
\begin{equation}\label{rho-b}
 \rho_n\leq1-r_nt_n
 =\rho_n+(1-\rho_n)(1-\Theta_n)
 \leq C_5\rho_n,
 \quad C_5:=1+2C_2+2C_4.
\end{equation}
Choosing $c_6=\frac{1}{C_5}$, we complete the proof of the lemma.
\end{proof}

\subsection{Moving-interval $L^p$-norm of $-\pa_xu_n$}
We now show that $\int^{\phi_n(t,q_n)}_{\phi_n(t,-q_n)}(-\partial_xu_n)^p(t,x)\dd x$ is sufficiently large, since $-\partial_x u_n(t,x)$ is uniformly large on $[\phi_n(t,-q_n),\phi_n(t,q_n)]$.
\begin{lemma}\label{yyh}
Let $n\gg 1$ and $p>2$. Denote $z_n(t,x)=-(\pa_xu_n)(t,x)$ and
\begin{equation*}
 \mathcal A_n(t)=\int_{I_n(t)}z_n^p(t,x)\dd x\quad \text{with}\quad \mathcal A_n(0)=2q_nn^p,
\end{equation*}
where the interval $I_n(t):=[\phi_n(t,-q_n),\phi_n(t,q_n)]$.
Then there exists a constant $\tilde{c}_0>0$ independent of $n$ such that
$$ \mathcal A_n(t_n)\geq \tilde{c}_0 q_nn^p\rho_n^{-(p-2)}.$$
\end{lemma}
\begin{proof}
In Eulerian variables, the second equation in \eqref{eq-u-ux} reads
\begin{align*}
 \pa_tz_n+u_n^Q \pa_xz_n=\frac Q2u^{Q-1}_nz_n^2-R_Q(u_n),
\end{align*}
which implies that for $p>2$, to be fixed later,
$$
\pa_t\left(z_n^p\right)+\pa_x\left(u_n^Qz_n^p\right)-\frac {Q(p-2)}2u^{Q-1}_nz_n^{p+1}=-pR_Q(u_n)z_n^{p-1}.
$$
Integrating the above over the interval $I_n(t)$ yields
\begin{align}\label{In(t)}
\int_{I_n(t)}\pa_t\left(z_n^p\right)\dd x&+\int_{I_n(t)}\pa_x\left(u_n^Qz_n^p\right)\dd x-\frac {Q(p-2)}2\int_{I_n(t)}u^{Q-1}_nz_n^{p+1}\dd x \nonumber
\\& =-p\int_{I_n(t)}R_Q(u_n)z_n^{p-1}\dd x.
\end{align}
Direct computations yield
\begin{align}\label{In(t)1}
\int_{I_n(t)}\pa_t\left(z_n^p\right)\dd x&=\frac{\dd}{\dd t}\int_{I_n(t)}z_n^p\dd x- u_n^{Q}(t,\phi(t,q_n))z_n^p(t,\phi(t,q_n))  \nonumber
\\& \qquad + u_n^Q(t,\phi(t,-q_n))z_n^p(t,\phi(t,-q_n)),
\end{align}
and
\begin{align}\label{In(t)2}
\int_{I_n(t)}\pa_x\left(u_n^Qz_n^p\right)\dd x &= u_n^Q(t,\phi(t,q_n))z_n^p(t,\phi(t,q_n)) \nonumber
\\& \qquad - u_n^Q(t,\phi(t,-q_n))z_n^p(t,\phi(t,-q_n)).
\end{align}
Substituting \eqref{In(t)1} and \eqref{In(t)2} into \eqref{In(t)}, we have
\begin{align}\label{An(t)-di=0}
 \mathcal A_n'(t)=(p-2)\int_{I_n(t)}\frac Q2u^{Q-1}_nz_n^{p+1}\dd x
       -p\int_{I_n(t)}R_Q(u_n)z_n^{p-1}\dd x.
\end{align}
By Corollary \ref{co-unQ-1}, for $0\leq t\leq t_n$,
\begin{align}\label{An(t)-di=1}
 \int_{I_n(t)}\frac Q2u^{Q-1}_nz_n^{p+1}\dd x
 \geq a_n^-\frac{n-d_n^-}{1-r_nt}\mathcal A_n(t).
\end{align}
It follows from \eqref{zn(t)-di} that $-(\pa_xu_{n})(t,x)
\geq n$ uniformly for $x\in[\phi(t,-q_n),\phi(t,q_n)]$, which yields
\begin{align}\label{An(t)-di=2}
 \left|\int_{I_n(t)}R_Q(u_n)z_n^{p-1}\dd x\right|
 \leq K_n\int_{I_n(t)}z_n^{p-1}\dd x
 \leq\frac{K_n}{n}\mathcal A_n(t).
\end{align}
Combining \eqref{An(t)-di=0}--\eqref{An(t)-di=2}, we obtain for $0\leq t\leq t_n$,
\begin{equation}\label{An(t)-di}
 \mathcal A_n'(t)
 \geq\left(\frac{(p-2)r_n}{1-r_nt}-\frac{pK_n}{n}\right)
           \mathcal A_n(t).
\end{equation}
Integrating this differential inequality yields the lower bound
\begin{equation}\label{eq:integral-growth}
 \mathcal A_n(t)
 \geq \mathcal A_n(0)\exp\!\left(-\frac{pK_nt}{n}\right)
                      (1-r_nt)^{-(p-2)},\qquad 0\leq t\leq t_n.
\end{equation}
Note that
\begin{align*}
 \frac{K_nt_n}{n}
 \leq\frac{K_n}{a_nn^2}
 =\frac{2C_0}{Q}\rho_n^2\longrightarrow0.
\end{align*}
and $\mathcal A_n(0)=2 q_nn^p$. Then from \eqref{rho-b}, we obtain
\begin{equation}\label{An(t)-geq}
 \mathcal A_n(t_n)\geq \tilde{c}_0 q_nn^p\rho_n^{-(p-2)},
\end{equation}
where $\tilde{c}_0>0$ is independent of $n$. This completes the proof of the lemma.
\end{proof}

\subsection{Norm inflation}
We now  prove that solutions to the generalized Camassa-Holm equation \eqref{gch} exhibits $H^s$-norm inflation on the timescale $t_n \approx \frac1{a_nn}
 =\frac2Qn^{-1+\kappa(Q-1)}$. It follows from Lemma \ref{le001} that
\begin{equation}\label{int-em}
\|u_{0,n}\|_{H^s}\leq \frac{C}{\log n}.
\end{equation}
Taking
$p=p_s$
in Lemma \ref{yyh}, applying \eqref{An(t)-geq} to $u_n$ and using the one-dimensional Sobolev embedding
$\norm{f_x}{L^{p_s}}\leq C_s\norm{f}{H^s}$ with $s-1\in(0,1/2)$ and $1/{p_s}=1/2-(s-1)$, we deduce that
\begin{align}
 \norm{u_n(t_n)}{H^s}
 &\geq C_s^{-1}\norm{\partial_xu_{n}(t_n)}{L^{p_s}}
 \geq c\,\mathcal A_n^{1/p_s}(t_n)\notag\\
 &\geq c\,nq_n^{1/p_s}\rho_n^{-(1-2/p_s)}\notag\\
 &=c\,\frac{n^{(1+\kappa)(1-2/p_s)}}{\log n}
 \longrightarrow\infty.
\end{align}
In conclusion, we obtain  norm inflation and hence  ill-posedness of the generalized Camassa-Holm equation \eqref{gch}. Thus, Theorem \ref{th2} is established.

\section*{Declarations}
\noindent\textbf{Data Availability}\\
No data was used for the research described in the article.

\vspace*{1em}
\noindent\textbf{Conflict of interest}\\
The authors declare that they have no conflict of interest.
\vspace*{1em}

\section*{Acknowledgements}

J. Li was supported by the National Natural Science Foundation of China (12661048), Ganpo Talents Project of Jiangxi Province (gpyc20240069) and Natural Science Foundation of Jiangxi Province (20252BAC210004 and 20252BAC240186). Xing Wu is supported by the Natural Science Foundation of Henan Province (252300421774). S. Zhang was supported by the Science and Technology Research Project of Jiangxi Provincial Department of Education(GJJ2401103) and Special Project for the Training of Early-Career Young Scientific and Technological Talents in Jiangxi Province (20252BEJ730255).


\end{document}